
%

\magnification=\magstep1
\catcode`\@=11
\font\tensmc=cmcsc10      
\def\smc{\tensmc}
\def\pagewidth#1{\hsize= #1 }
\def\pageheight#1{\vsize= #1 }
\def\hcorrection#1{\advance\hoffset by #1 }
\def\vcorrection#1{\advance\voffset by #1 }
\def\wlog#1{}
\newif\iftitle@
\outer\def\title{\title@true\vglue 24\p@ plus 12\p@ minus 12\p@
   \bgroup\let\\=\cr\tabskip\centering
   \halign to \hsize\bgroup\tenbf\hfill\ignorespaces##\unskip\hfill\cr}
\def\endtitle{\cr\egroup\egroup\vglue 18\p@ plus 12\p@ minus 6\p@}
\outer\def\author{\iftitle@\vglue -18\p@ plus -12\p@ minus -6\p@\fi\vglue
    12\p@ plus 6\p@ minus 3\p@\bgroup\let\\=\cr\tabskip\centering
    \halign to \hsize\bgroup\smc\hfill\ignorespaces##\unskip\hfill\cr}
\def\endauthor{\cr\egroup\egroup\vglue 18\p@ plus 12\p@ minus 6\p@}
\outer\def\heading{\bigbreak\bgroup\let\\=\cr\tabskip\centering
    \halign to \hsize\bgroup\smc\hfill\ignorespaces##\unskip\hfill\cr}
\def\endheading{\cr\egroup\egroup\nobreak\medskip}

\outer\def\endproclaim{\par\ifdim\lastskip<\medskipamount\removelastskip
  \penalty 55 \fi\medskip\rm}
\outer\def\demo#1{\par\ifdim\lastskip<\smallskipamount\removelastskip
    \smallskip\fi\noindent{\smc\ignorespaces#1\unskip:\enspace}\rm
      \ignorespaces}

\newcount\footmarkcount@
\footmarkcount@=1
\def\makefootnote@#1#2{\insert\footins{\interlinepenalty=100
  \splittopskip=\ht\strutbox \splitmaxdepth=\dp\strutbox
  \floatingpenalty=\@MM
  \leftskip=\z@\rightskip=\z@\spaceskip=\z@\xspaceskip=\z@
  \noindent{#1}\footstrut\rm\ignorespaces #2\strut}}
\def\footnote{\let\@sf=\empty\ifhmode\edef\@sf{\spacefactor
   =\the\spacefactor}\/\fi\futurelet\next\footnote@}
\def\footnote@{\ifx"\next\let\next\footnote@@\else
    \let\next\footnote@@@\fi\next}
\def\footnote@@"#1"#2{#1\@sf\relax\makefootnote@{#1}{#2}}
\def\footnote@@@#1{$^{\number\footmarkcount@}$\makefootnote@
   {$^{\number\footmarkcount@}$}{#1}\global\advance\footmarkcount@ by 1 }

\hyphenation{man-u-script man-u-scripts ap-pen-dix ap-pen-di-ces}
\hyphenation{data-base data-bases}
\ifx\amstexloaded@\relax\catcode`\@=13
  \endinput\else\let\amstexloaded@=\relax\fi
\newlinechar=`\^^J
\def\eat@#1{}
\def\Space@.{\futurelet\Space@\relax}
\Space@. %
\newhelp\athelp@
{Only certain combinations beginning with @ make sense to me.^^J
Perhaps you wanted \string\@\space for a printed @?^^J
I've ignored the character or group after @.}
\def\futureletnextat@{\futurelet\next\at@}
{\catcode`\@=\active
\lccode`\Z=`\@ \lowercase
{\gdef@{\expandafter\csname futureletnextatZ\endcsname}
\expandafter\gdef\csname atZ\endcsname
   {\ifcat\noexpand\next a\def\next{\csname atZZ\endcsname}\else
   \ifcat\noexpand\next0\def\next{\csname atZZ\endcsname}\else
    \def\next{\csname atZZZ\endcsname}\fi\fi\next}
\expandafter\gdef\csname atZZ\endcsname#1{\expandafter
   \ifx\csname #1Zat\endcsname\relax\def\next
     {\errhelp\expandafter=\csname athelpZ\endcsname
      \errmessage{Invalid use of \string@}}\else
       \def\next{\csname #1Zat\endcsname}\fi\next}
\expandafter\gdef\csname atZZZ\endcsname#1{\errhelp
    \expandafter=\csname athelpZ\endcsname
      \errmessage{Invalid use of \string@}}}}
\def\atdef@#1{\expandafter\def\csname #1@at\endcsname}
\newhelp\defahelp@{If you typed \string\define\space cs instead of
\string\define\string\cs\space^^J
I've substituted an inaccessible control sequence so that your^^J
definition will be completed without mixing me up too badly.^^J
If you typed \string\define{\string\cs} the inaccessible control sequence^^J
was defined to be \string\cs, and the rest of your^^J
definition appears as input.}
\newhelp\defbhelp@{I've ignored your definition, because it might^^J
conflict with other uses that are important to me.}
\def\define{\futurelet\next\define@}
\def\define@{\ifcat\noexpand\next\relax
  \def\next{\define@@}%
  \else\errhelp=\defahelp@
  \errmessage{\string\define\space must be followed by a control
     sequence}\def\next{\def\garbage@}\fi\next}
\def\undefined@{}
\def\preloaded@{}
\def\define@@#1{\ifx#1\relax\errhelp=\defbhelp@
   \errmessage{\string#1\space is already defined}\def\next{\def\garbage@}%
   \else\expandafter\ifx\csname\expandafter\eat@\string
         #1@\endcsname\undefined@\errhelp=\defbhelp@
   \errmessage{\string#1\space can't be defined}\def\next{\def\garbage@}%
   \else\expandafter\ifx\csname\expandafter\eat@\string#1\endcsname\relax
     \def\next{\def#1}\else\errhelp=\defbhelp@
     \errmessage{\string#1\space is already defined}\def\next{\def\garbage@}%
      \fi\fi\fi\next}
\def\famzero{\fam\z@}

\def\det{\mathop{\famzero det}}
\def\dim{\mathop{\famzero dim}\nolimits}
\def\exp{\mathop{\famzero exp}\nolimits}

\def\inf{\mathop{\famzero inf}}

\def\lim{\mathop{\famzero lim}}

\def\log{\mathop{\famzero log}\nolimits}

\def\min{\mathop{\famzero min}}

\def\sup{\mathop{\famzero sup}}

\def\textfont@#1#2{\def#1{\relax\ifmmode
    \errmessage{Use \string#1\space only in text}\else#2\fi}}
\textfont@\rm\tenrm
\textfont@\it\tenit
\textfont@\sl\tensl
\textfont@\bf\tenbf
\textfont@\smc\tensmc
\let\ic@=\/
\def\/{\unskip\ic@}
\def\textfonti{\the\textfont1 }
\def\t#1#2{{\edef\next{\the\font}\textfonti\accent"7F \next#1#2}}
\let\B=\=
\let\D=\.
\def~{\unskip\nobreak\ \ignorespaces}
{\catcode`\@=\active
\gdef\@{\char'100 }}
\atdef@-{\leavevmode\futurelet\next\athyph@}
\def\athyph@{\ifx\next-\let\next=\athyph@@
  \else\let\next=\athyph@@@\fi\next}
\def\athyph@@@{\hbox{-}}
\def\athyph@@#1{\futurelet\next\athyph@@@@}
\def\athyph@@@@{\if\next-\def\next##1{\hbox{---}}\else
    \def\next{\hbox{--}}\fi\next}
\def\.{.\spacefactor=\@m}
\atdef@.{\null.}
\atdef@,{\null,}
\atdef@;{\null;}
\atdef@:{\null:}
\atdef@?{\null?}
\atdef@!{\null!}
\def\srdr@{\thinspace}
\def\drsr@{\kern.02778em}
\def\sldl@{\kern.02778em}
\def\dlsl@{\thinspace}
\atdef@"{\unskip\futurelet\next\atqq@}
\def\atqq@{\ifx\next\Space@\def\next. {\atqq@@}\else
         \def\next.{\atqq@@}\fi\next.}
\def\atqq@@{\futurelet\next\atqq@@@}
\def\atqq@@@{\ifx\next`\def\next`{\atqql@}\else\def\next'{\atqqr@}\fi\next}
\def\atqql@{\futurelet\next\atqql@@}
\def\atqql@@{\ifx\next`\def\next`{\sldl@``}\else\def\next{\dlsl@`}\fi\next}
\def\atqqr@{\futurelet\next\atqqr@@}
\def\atqqr@@{\ifx\next'\def\next'{\srdr@''}\else\def\next{\drsr@'}\fi\next}

\def\textfontii{\the\textfont2 }
\def\{{\relax\ifmmode\lbrace\else
    {\textfontii f}\spacefactor=\@m\fi}
\def\}{\relax\ifmmode\rbrace\else
    \let\@sf=\empty\ifhmode\edef\@sf{\spacefactor=\the\spacefactor}\fi
      {\textfontii g}\@sf\relax\fi}
\def\nonhmodeerr@#1{\errmessage
     {\string#1\space allowed only within text}}
\def\linebreak{\relax\ifhmode\unskip\break\else
    \nonhmodeerr@\linebreak\fi}
\def\allowlinebreak{\relax
   \ifhmode\allowbreak\else\nonhmodeerr@\allowlinebreak\fi}
\newskip\saveskip@
\def\nolinebreak{\relax\ifhmode\saveskip@=\lastskip\unskip
  \nobreak\ifdim\saveskip@>\z@\hskip\saveskip@\fi
   \else\nonhmodeerr@\nolinebreak\fi}
\def\newline{\relax\ifhmode\null\hfil\break
    \else\nonhmodeerr@\newline\fi}
\def\nonmathaerr@#1{\errmessage
     {\string#1\space is not allowed in display math mode}}
\def\nonmathberr@#1{\errmessage{\string#1\space is allowed only in math mode}}
\def\mathbreak{\relax\ifmmode\ifinner\break\else
   \nonmathaerr@\mathbreak\fi\else\nonmathberr@\mathbreak\fi}
\def\nomathbreak{\relax\ifmmode\ifinner\nobreak\else
    \nonmathaerr@\nomathbreak\fi\else\nonmathberr@\nomathbreak\fi}
\def\allowmathbreak{\relax\ifmmode\ifinner\allowbreak\else
     \nonmathaerr@\allowmathbreak\fi\else\nonmathberr@\allowmathbreak\fi}
\def\pagebreak{\relax\ifmmode
   \ifinner\errmessage{\string\pagebreak\space
     not allowed in non-display math mode}\else\postdisplaypenalty-\@M\fi
   \else\ifvmode\penalty-\@M\else\edef\spacefactor@
       {\spacefactor=\the\spacefactor}\vadjust{\penalty-\@M}\spacefactor@
        \relax\fi\fi}
\def\nopagebreak{\relax\ifmmode
     \ifinner\errmessage{\string\nopagebreak\space
    not allowed in non-display math mode}\else\postdisplaypenalty\@M\fi
    \else\ifvmode\nobreak\else\edef\spacefactor@
        {\spacefactor=\the\spacefactor}\vadjust{\penalty\@M}\spacefactor@
         \relax\fi\fi}
\def\newpage{\relax\ifvmode\vfill\penalty-\@M\else\nonvmodeerr@\newpage\fi}
\def\nonvmodeerr@#1{\errmessage
    {\string#1\space is allowed only between paragraphs}}
\def\smallpagebreak{\relax\ifvmode\smallbreak
      \else\nonvmodeerr@\smallpagebreak\fi}
\def\medpagebreak{\relax\ifvmode\medbreak
       \else\nonvmodeerr@\medpagebreak\fi}
\def\bigpagebreak{\relax\ifvmode\bigbreak
      \else\nonvmodeerr@\bigpagebreak\fi}
\newdimen\captionwidth@
\captionwidth@=\hsize
\advance\captionwidth@ by -1.5in
\def\caption#1{}
\def\topspace#1{\gdef\thespace@{#1}\ifvmode\def\next
    {\futurelet\next\topspace@}\else\def\next{\nonvmodeerr@\topspace}\fi\next}
\def\topspace@{\ifx\next\Space@\def\next. {\futurelet\next\topspace@@}\else
     \def\next.{\futurelet\next\topspace@@}\fi\next.}
\def\topspace@@{\ifx\next\caption\let\next\topspace@@@\else
    \let\next\topspace@@@@\fi\next}
 \def\topspace@@@@{\topinsert\vbox to
       \thespace@{}\endinsert}
\def\topspace@@@\caption#1{\topinsert\vbox to
    \thespace@{}\nobreak
      \smallskip
    \setbox\z@=\hbox{\noindent\ignorespaces#1\unskip}%
   \ifdim\wd\z@>\captionwidth@
   \centerline{\vbox{\hsize=\captionwidth@\noindent\ignorespaces#1\unskip}}%
   \else\centerline{\box\z@}\fi\endinsert}
\def\midspace#1{\gdef\thespace@{#1}\ifvmode\def\next
    {\futurelet\next\midspace@}\else\def\next{\nonvmodeerr@\midspace}\fi\next}
\def\midspace@{\ifx\next\Space@\def\next. {\futurelet\next\midspace@@}\else
     \def\next.{\futurelet\next\midspace@@}\fi\next.}
\def\midspace@@{\ifx\next\caption\let\next\midspace@@@\else
    \let\next\midspace@@@@\fi\next}
 \def\midspace@@@@{\midinsert\vbox to
       \thespace@{}\endinsert}
\def\midspace@@@\caption#1{\midinsert\vbox to
    \thespace@{}\nobreak
      \smallskip
      \setbox\z@=\hbox{\noindent\ignorespaces#1\unskip}%
      \ifdim\wd\z@>\captionwidth@
    \centerline{\vbox{\hsize=\captionwidth@\noindent\ignorespaces#1\unskip}}%
    \else\centerline{\box\z@}\fi\endinsert}
\mathchardef\prime@="0230
\def\prime{{{}\prime@{}}}
\def\prim@s{\prime@\futurelet\next\pr@m@s}

\def\,{\relax\ifmmode\mskip\thinmuskip\else\thinspace\fi}
\def\!{\relax\ifmmode\mskip-\thinmuskip\else\negthinspace\fi}
\def\frac#1#2{{#1\over#2}}

\def\:{\nobreak\hskip.1111em{:}\hskip.3333em plus .0555em\relax}
\def\intic@{\mathchoice{\hskip5\p@}{\hskip4\p@}{\hskip4\p@}{\hskip4\p@}}
\def\negintic@
 {\mathchoice{\hskip-5\p@}{\hskip-4\p@}{\hskip-4\p@}{\hskip-4\p@}}
\def\intkern@{\mathchoice{\!\!\!}{\!\!}{\!\!}{\!\!}}
\def\intdots@{\mathchoice{\cdots}{{\cdotp}\mkern1.5mu
    {\cdotp}\mkern1.5mu{\cdotp}}{{\cdotp}\mkern1mu{\cdotp}\mkern1mu
      {\cdotp}}{{\cdotp}\mkern1mu{\cdotp}\mkern1mu{\cdotp}}}
\newcount\intno@
\def\iint{\intno@=\tw@\futurelet\next\ints@}
\def\iiint{\intno@=\thr@@\futurelet\next\ints@}
\def\iiiint{\intno@=4 \futurelet\next\ints@}
\def\idotsint{\intno@=\z@\futurelet\next\ints@}
\def\ints@{\findlimits@\ints@@}
\newif\iflimtoken@
\newif\iflimits@
\def\findlimits@{\limtoken@false\limits@false\ifx\next\limits
 \limtoken@true\limits@true\else\ifx\next\nolimits\limtoken@true\limits@false
    \fi\fi}
\def\multintlimits@{\intop\ifnum\intno@=\z@\intdots@
  \else\intkern@\fi
    \ifnum\intno@>\tw@\intop\intkern@\fi
     \ifnum\intno@>\thr@@\intop\intkern@\fi\intop}
\def\multint@{\int\ifnum\intno@=\z@\intdots@\else\intkern@\fi
   \ifnum\intno@>\tw@\int\intkern@\fi
    \ifnum\intno@>\thr@@\int\intkern@\fi\int}
\def\ints@@{\iflimtoken@\def\ints@@@{\iflimits@
   \negintic@\mathop{\intic@\multintlimits@}\limits\else
    \multint@\nolimits\fi\eat@}\else
     \def\ints@@@{\multint@\nolimits}\fi\ints@@@}
\def\Sb{_\bgroup\vspace@
        \baselineskip=\fontdimen10 \scriptfont\tw@
        \advance\baselineskip by \fontdimen12 \scriptfont\tw@
        \lineskip=\thr@@\fontdimen8 \scriptfont\thr@@
        \lineskiplimit=\thr@@\fontdimen8 \scriptfont\thr@@
        \Let@\vbox\bgroup\halign\bgroup \hfil$\scriptstyle
            {##}$\hfil\cr}
\def\endSb{\crcr\egroup\egroup\egroup}
\def\Sp{^\bgroup\vspace@
        \baselineskip=\fontdimen10 \scriptfont\tw@
        \advance\baselineskip by \fontdimen12 \scriptfont\tw@
        \lineskip=\thr@@\fontdimen8 \scriptfont\thr@@
        \lineskiplimit=\thr@@\fontdimen8 \scriptfont\thr@@
        \Let@\vbox\bgroup\halign\bgroup \hfil$\scriptstyle
            {##}$\hfil\cr}
\def\endSp{\crcr\egroup\egroup\egroup}
\def\Let@{\relax\iffalse{\fi\let\\=\cr\iffalse}\fi}
\def\vspace@{\def\vspace##1{\noalign{\vskip##1 }}}
\def\aligned{\,\vcenter\bgroup\vspace@\Let@\openup\jot\m@th\ialign
  \bgroup \strut\hfil$\displaystyle{##}$&$\displaystyle{{}##}$\hfil\crcr}
\def\endaligned{\crcr\egroup\egroup}
\def\matrix{\,\vcenter\bgroup\Let@\vspace@
    \normalbaselines
  \m@th\ialign\bgroup\hfil$##$\hfil&&\quad\hfil$##$\hfil\crcr
    \mathstrut\crcr\noalign{\kern-\baselineskip}}
\def\endmatrix{\crcr\mathstrut\crcr\noalign{\kern-\baselineskip}\egroup
                \egroup\,}
\newtoks\hashtoks@
\hashtoks@={#}
\def\format{\crcr\egroup\iffalse{\fi\ifnum`}=0 \fi\format@}
\def\format@#1\\{\def\preamble@{#1}%
  \def\c{\hfil$\the\hashtoks@$\hfil}%
  \def\r{\hfil$\the\hashtoks@$}%
  \def\l{$\the\hashtoks@$\hfil}%
  \setbox\z@=\hbox{\xdef\Preamble@{\preamble@}}\ifnum`{=0 \fi\iffalse}\fi
   \ialign\bgroup\span\Preamble@\crcr}
\def\pmatrix{\left(\matrix} \def\endpmatrix{\endmatrix\right)}

\def\cases{\left\{\,\vcenter\bgroup\vspace@
     \normalbaselines\openup\jot\m@th
       \Let@\ialign\bgroup$##$\hfil&\quad$##$\hfil\crcr
      \mathstrut\crcr\noalign{\kern-\baselineskip}}

\newif\iftagsleft@
\tagsleft@true
\def\TagsOnRight{\global\tagsleft@false}
\def\tag#1$${\iftagsleft@\leqno\else\eqno\fi
 \hbox{\def\pagebreak{\global\postdisplaypenalty-\@M}%
 \def\nopagebreak{\global\postdisplaypenalty\@M}\rm(#1\unskip)}%
  $$\postdisplaypenalty\z@\ignorespaces}
\interdisplaylinepenalty=\@M
\def\allowdisplaybreak@{\def\allowdisplaybreak{\noalign{\allowbreak}}}
\def\displaybreak@{\def\displaybreak{\noalign{\break}}}
\def\align#1\endalign{\def\tag{&}\vspace@\allowdisplaybreak@\displaybreak@
  \iftagsleft@\lalign@#1\endalign\else
   \ralign@#1\endalign\fi}
\def\ralign@#1\endalign{\displ@y\Let@\tabskip\centering\halign to\displaywidth
     {\hfil$\displaystyle{##}$\tabskip=\z@&$\displaystyle{{}##}$\hfil
       \tabskip=\centering&\llap{\hbox{(\rm##\unskip)}}\tabskip\z@\crcr
             #1\crcr}}
\def\lalign@
 #1\endalign{\displ@y\Let@\tabskip\centering\halign to \displaywidth
   {\hfil$\displaystyle{##}$\tabskip=\z@&$\displaystyle{{}##}$\hfil
   \tabskip=\centering&\kern-\displaywidth
        \rlap{\hbox{(\rm##\unskip)}}\tabskip=\displaywidth\crcr
               #1\crcr}}
\def\overrightarrow{\mathpalette\overrightarrow@}
\def\overrightarrow@#1#2{\vbox{\ialign{$##$\cr
    #1{-}\mkern-6mu\cleaders\hbox{$#1\mkern-2mu{-}\mkern-2mu$}\hfill
     \mkern-6mu{\to}\cr
     \noalign{\kern -1\p@\nointerlineskip}
     \hfil#1#2\hfil\cr}}}
\def\overleftarrow{\mathpalette\overleftarrow@}
\def\overleftarrow@#1#2{\vbox{\ialign{$##$\cr
     #1{\leftarrow}\mkern-6mu\cleaders\hbox{$#1\mkern-2mu{-}\mkern-2mu$}\hfill
      \mkern-6mu{-}\cr
     \noalign{\kern -1\p@\nointerlineskip}
     \hfil#1#2\hfil\cr}}}
\def\overleftrightarrow{\mathpalette\overleftrightarrow@}
\def\overleftrightarrow@#1#2{\vbox{\ialign{$##$\cr
     #1{\leftarrow}\mkern-6mu\cleaders\hbox{$#1\mkern-2mu{-}\mkern-2mu$}\hfill
       \mkern-6mu{\to}\cr
    \noalign{\kern -1\p@\nointerlineskip}
      \hfil#1#2\hfil\cr}}}
\def\underrightarrow{\mathpalette\underrightarrow@}
\def\underrightarrow@#1#2{\vtop{\ialign{$##$\cr
    \hfil#1#2\hfil\cr
     \noalign{\kern -1\p@\nointerlineskip}
    #1{-}\mkern-6mu\cleaders\hbox{$#1\mkern-2mu{-}\mkern-2mu$}\hfill
     \mkern-6mu{\to}\cr}}}
\def\underleftarrow{\mathpalette\underleftarrow@}
\def\underleftarrow@#1#2{\vtop{\ialign{$##$\cr
     \hfil#1#2\hfil\cr
     \noalign{\kern -1\p@\nointerlineskip}
     #1{\leftarrow}\mkern-6mu\cleaders\hbox{$#1\mkern-2mu{-}\mkern-2mu$}\hfill
      \mkern-6mu{-}\cr}}}
\def\underleftrightarrow{\mathpalette\underleftrightarrow@}
\def\underleftrightarrow@#1#2{\vtop{\ialign{$##$\cr
      \hfil#1#2\hfil\cr
    \noalign{\kern -1\p@\nointerlineskip}
     #1{\leftarrow}\mkern-6mu\cleaders\hbox{$#1\mkern-2mu{-}\mkern-2mu$}\hfill
       \mkern-6mu{\to}\cr}}}
\def\sqrt#1{\radical"270370 {#1}}
\def\dots{\relax\ifmmode\let\next=\ldots\else\let\next=\tdots@\fi\next}
\def\tdots@{\unskip\ \tdots@@}
\def\tdots@@{\futurelet\next\tdots@@@}
\def\tdots@@@{$\mathinner{\ldotp\ldotp\ldotp}\,
   \ifx\next,$\else
   \ifx\next.\,$\else
   \ifx\next;\,$\else
   \ifx\next:\,$\else
   \ifx\next?\,$\else
   \ifx\next!\,$\else
   $ \fi\fi\fi\fi\fi\fi}
\def\text{\relax\ifmmode\let\next=\text@\else\let\next=\text@@\fi\next}
\def\text@@#1{\hbox{#1}}
\def\text@#1{\mathchoice
 {\hbox{\everymath{\displaystyle}\def\textfonti{\the\textfont1 }%
    \def\textfontii{\the\textfont2 }\textdef@@ T#1}}
 {\hbox{\everymath{\textstyle}\def\textfonti{\the\textfont1 }%
    \def\textfontii{\the\textfont2 }\textdef@@ T#1}}
 {\hbox{\everymath{\scriptstyle}\def\textfonti{\the\scriptfont1 }%
   \def\textfontii{\the\scriptfont2 }\textdef@@ S\rm#1}}
 {\hbox{\everymath{\scriptscriptstyle}\def\textfonti{\the\scriptscriptfont1 }%
   \def\textfontii{\the\scriptscriptfont2 }\textdef@@ s\rm#1}}}
\def\textdef@@#1{\textdef@#1\rm \textdef@#1\bf
   \textdef@#1\sl \textdef@#1\it}

\def\textdef@#1#2{\def\next{\csname\expandafter\eat@\string#2fam\endcsname}%
\if S#1\edef#2{\the\scriptfont\next\relax}%
 \else\if s#1\edef#2{\the\scriptscriptfont\next\relax}%
 \else\edef#2{\the\textfont\next\relax}\fi\fi}
\scriptfont\itfam=\tenit \scriptscriptfont\itfam=\tenit
\scriptfont\slfam=\tensl \scriptscriptfont\slfam=\tensl
\mathcode`\0="0030
\mathcode`\1="0031
\mathcode`\2="0032
\mathcode`\3="0033
\mathcode`\4="0034
\mathcode`\5="0035
\mathcode`\6="0036
\mathcode`\7="0037
\mathcode`\8="0038
\mathcode`\9="0039
\def\Cal{\relax\ifmmode\let\next=\Cal@\else
     \def\next{\errmessage{Use \string\Cal\space only in math mode}}\fi\next}
\def\Cal@#1{{\fam2 #1}}
\def\bold{\relax\ifmmode\let\next=\bold@\else
   \def\next{\errmessage{Use \string\bold\space only in math
      mode}}\fi\next}\def\bold@#1{{\fam\bffam #1}}
\mathchardef\Gamma="0000
\mathchardef\Delta="0001
\mathchardef\Theta="0002
\mathchardef\Lambda="0003
\mathchardef\Xi="0004
\mathchardef\Pi="0005
\mathchardef\Sigma="0006
\mathchardef\Upsilon="0007
\mathchardef\Phi="0008
\mathchardef\Psi="0009
\mathchardef\Omega="000A
\mathchardef\varGamma="0100
\mathchardef\varDelta="0101
\mathchardef\varTheta="0102
\mathchardef\varLambda="0103
\mathchardef\varXi="0104
\mathchardef\varPi="0105
\mathchardef\varSigma="0106
\mathchardef\varUpsilon="0107
\mathchardef\varPhi="0108
\mathchardef\varPsi="0109
\mathchardef\varOmega="010A
\font\dummyft@=dummy
\fontdimen1 \dummyft@=\z@
\fontdimen2 \dummyft@=\z@
\fontdimen3 \dummyft@=\z@
\fontdimen4 \dummyft@=\z@
\fontdimen5 \dummyft@=\z@
\fontdimen6 \dummyft@=\z@
\fontdimen7 \dummyft@=\z@
\fontdimen8 \dummyft@=\z@
\fontdimen9 \dummyft@=\z@
\fontdimen10 \dummyft@=\z@
\fontdimen11 \dummyft@=\z@
\fontdimen12 \dummyft@=\z@
\fontdimen13 \dummyft@=\z@
\fontdimen14 \dummyft@=\z@
\fontdimen15 \dummyft@=\z@
\fontdimen16 \dummyft@=\z@
\fontdimen17 \dummyft@=\z@
\fontdimen18 \dummyft@=\z@
\fontdimen19 \dummyft@=\z@
\fontdimen20 \dummyft@=\z@
\fontdimen21 \dummyft@=\z@
\fontdimen22 \dummyft@=\z@
\def\fontlist@{\\{\tenrm}\\{\sevenrm}\\{\fiverm}\\{\teni}\\{\seveni}%
 \\{\fivei}\\{\tensy}\\{\sevensy}\\{\fivesy}\\{\tenex}\\{\tenbf}\\{\sevenbf}%
 \\{\fivebf}\\{\tensl}\\{\tenit}\\{\tensmc}}
\def\dodummy@{{\def\\##1{\global\let##1=\dummyft@}\fontlist@}}
\newif\ifsyntax@
\newcount\countxviii@
\def\newtoks@{\alloc@5\toks\toksdef\@cclvi}
\def\nopages@{\output={\setbox\z@=\box\@cclv \deadcycles=\z@}\newtoks@\output}
\def\syntax{\syntax@true\dodummy@\countxviii@=\count18
\loop \ifnum\countxviii@ > \z@ \textfont\countxviii@=\dummyft@
   \scriptfont\countxviii@=\dummyft@ \scriptscriptfont\countxviii@=\dummyft@
     \advance\countxviii@ by-\@ne\repeat
\dummyft@\tracinglostchars=\z@
  \nopages@\frenchspacing\hbadness=\@M}
\def\magstep#1{\ifcase#1 1000\or
 1200\or 1440\or 1728\or 2074\or 2488\or
 \errmessage{\string\magstep\space only works up to 5}\fi\relax}
{\lccode`\2=`\p \lccode`\3=`\t
 \lowercase{\gdef\tru@#123{#1truept}}}

\def\scaletype#1{\mag=#1\relax
 \hsize=\expandafter\tru@\the\hsize
 \vsize=\expandafter\tru@\the\vsize
 \dimen\footins=\expandafter\tru@\the\dimen\footins}

\def\scalefont#1#2\andcallit#3{\edef\font@{\the\font}#1\font#3=
  \fontname\font\space scaled #2\relax\font@}
\def\Mag@#1#2{\ifdim#1<1pt\multiply#1 #2\relax\divide#1 1000 \else
  \ifdim#1<10pt\divide#1 10 \multiply#1 #2\relax\divide#1 100\else
  \divide#1 100 \multiply#1 #2\relax\divide#1 10 \fi\fi}
\def\scalelinespacing#1{\Mag@\baselineskip{#1}\Mag@\lineskip{#1}%
  \Mag@\lineskiplimit{#1}}
\def\wlog#1{\immediate\write-1{#1}}
\catcode`\@=\active



\baselineskip=12pt
\pagewidth{4.5 in}\pageheight{7.3 in}
\hoffset=.312 in \voffset=.437 in


\font\Lbf=cmbx10 scaled\magstep1
\font\bi=cmbxti10
\font\sc=cmcsc10
\font\Srm=cmr9


\def\a{\alpha}

\def\e{\varepsilon}
\def\f{\varphi}

\def\l{\lambda}
\def\p{\psi}
\def\r{\rho}

\def\t{\tau}
\def\dd{\Delta}


\def\hs{Hilbert space}
\def\op{operator}

\def\al{algebra}

\def\ca{\hbox{C*-\al}}

\def\hm{homomorphism}

\def\i{\infty}

\def\fd{finite-dimensional}

\def\el{element}
\def\cs{\hbox{C*-subalgebra}}

\def\uca{unital \ca}

\def\ifo{if and only if}


\def\tn{{\bold T}}
\def\cn{{\bold C}}
\def\zn{{\bold Z}}
\def\rn{{\bold R}}

\def\nn{{\bold N}}

\def\n#1{{\|#1\|}}
\def\m#1{{\vert #1 \vert}}

\def\ss{\subseteq}

\def\inv#1{{#1}^{-1}}

\def\oper#1{\mathop{\text{\rm #1}}}

\def\limn{\lim_{n\to \i}}

\def\det{\oper{det}}

\def\dim{\oper{dim}}

\def\Inv{\oper{Inv}}


\def\ct{\centerline}
\def\ni{\noindent}
\def\sm{\smallskip}
\def\md{\medskip}
\def\bg{\bigskip}


\def\tf{therefore}

\def\st{such that}
\def\tes{there exists}


\hyphenation{theorem theory hyper-finite tracial nuclear
pro-ject-ion pro-ject-ions sym-metry sym-metries gen-erate
gen-erates gen-erat-ed unit-al anal-yt-ic
oper-ator oper-ators iso-morph-ism homo-morph-ism acad-em-ic
di-vis-or di-vis-ors top-olog-ic-al suprem-um stable equi-val-ent}


\def\smalltitle{\lowercase\expandafter{\Title}}
\headline={\ifnum\pageno=1 \hfill
\else \vbox{\line{\sc\smalltitle\hfill{\bf\folio}}%
\vskip 1.5pt\hrule}\hss\fi}

\footline={\ifnum\pageno=1 \hfill\folio\hfill \else \hfill\fi}

\hfuzz=1pt

\def\maketitle{\null\par\bg\bg
\ct{\bf \Title}
\ct{\sc G.J. Murphy and N.C. Phillips\footnote{Partially supported by
NSF grant \grant \newline AMS subject classification: \ams }}
\bg\bg}

\def\makelongtitle{\null\par\bg\bg
\ct{\Lbf \lineone}
\ct{\Lbf \linetwo}\bg
\ct{\sc G.J. Murphy and N.C. Phillips\footnote{Partially supported by
NSF grant \grant \newline AMS subject classification: \ams }}
\bg\bg}

\def\abstract#1{{\narrower\Srm
\advance\baselineskip by -2pt\ni{\bf Abstract\. }#1\par}\qedd}

\def\sect#1{\thmno=0\advance\scno by 1
\bigpagebreak\ni{\Lbf \the\scno.~#1}\md\nobreak}

\def\protothm#1{\advance\thmno by1
\sl\bigpagebreak\ni{\bf
\if\setsections Y\the\scno.\else \fi\the\thmno.~#1. }\,}

\def\thm{\protothm{Theorem}}
\def\lem{\protothm{Lemma}}

\def\pro{\protothm{Proposition}}

\def\exa{\advance\thmno by1
\bigpagebreak\ni{\bf
\if\setsections Y\the\scno.\else \fi\the\thmno.~{\bi Example.} }\,}
\def\rem{\advance\thmno by1
\bigpagebreak\ni{\bf
\if\setsections Y\the\scno.\else \fi\the\thmno.~{\bi Remark.} }\,}

\def\pf{\rm\medpagebreak\ni{\bi Proof. }\,}
\def\qed{~$\sqcup$\llap{$\sqcap$}\par\medpagebreak}
\def\qedd{\rm\par\medpagebreak}

\def\1{\nobreak\sm\item{\rm (1)}}
\def\2{\item{\rm (2)}}
\def\3{\item{\rm (3)}}
\def\4{\item{\rm (4)}}
\def\5{\item{\rm (5)}}
\def\6{\item{\rm (6)}}

\def\refs{\bigpagebreak\ni{\Lbf References}\md\nobreak
\hfuzz=3pt\overfullrule=0pt\frenchspacing}
\def\refbox{\ifnum\refno>9 \item{\phantom 9[\the\refno]}
           \else \item{[\the\refno]}\fi}
\def\pref#1#2#3#4#5#6{\advance\refno by1
\refbox #1, #2, {\it #3 }{\bf #4 }(19#5), #6.}
\def\ppref#1#2#3#4{\advance\refno by1
\refbox #1, #2, preprint (19#3) #4.}
\def\psub#1#2{\advance\refno by1
\refbox #1, #2, submitted.}
\def\bref#1#2#3#4#5{\advance \refno by 1
\refbox #1, {\sl #2. }#3, #4, 19#5.}
\def\thesis#1#2#3#4{\advance\refno by1
\refbox #1, #2 (PhD thesis),  #3, 19#4.}

\def\address{\nobreak\bg{\sc
\baselineskip=12pt\vbox{\line{Department of Mathematics, \hfill}
\line{University College,\hfill}
\line{Cork, Ireland.\hfill}
\line{\,\,\,\,\,\, and \hfill}
\line{Department of Mathematics, \hfill}
\line{University of Oregon, \hfill}
\line{Eugene OR 97403-1222, U.S.A. \hfill}}}}

\newcount\scno
\newcount\thmno
\newcount\refno
\scno=0
\thmno=0
\refno=0


\def\Title{The Approximate Positive Factorization Property}
\def\ams{46L05, 46L10}
\def\setsections{Y}
\def\grant{DMS-9106285.}



\def\lineone{C*-Algebras with the Approximate Positive}
\def\linetwo{Factorization Property}
\makelongtitle


\def\BAL{1} \def\Bll{\BAL}
\def\BDR{2}
\def\BP{3}
\def\CHO{4}
\def\Choi{5}
\def\Cu{6}
\def\CP{7}
\def\dlHS{8}
\def\GdB{9}
\def\Gd{10}
\def\Hg{11}
\def\HAL{12}
\def\Ks{13}
\def\KLMR{14}
\def\Leen{15}
\def\Ln{16}
\def\MUR{17}
\def\PT{18}
\def\Ph{19}
\def\PhB{20}
\def\PhC{21}
\def\QUI{22}
\def\RAD{23}
\def\RIE{24} \def\RfA{\RIE}
\def\Rf{25}
\def\RrA{26}
\def\RrB{27}
\def\RS{28}
\def\WU{29}

\def\Inv{\oper{Inv}}
\def\Invo{\Inv_0}
\def\rk{\oper{rank}}

\def\APFP{approximate positive factorization property}
\def\mops{mutually orthogonal projections}
\def\CPP#1{P( #1)^-}
\def\dirlim{\displaystyle \lim_{\longrightarrow}}
\def\zt{\zeta}

\def\dfn{\advance\thmno by1
\bigpagebreak\ni{\bf
\if\setsections Y\the\scno.\else \fi\the\thmno.~{\bi Definition.} }\,}
\def\ntn{\advance\thmno by1
\bigpagebreak\ni{\bf
\if\setsections Y\the\scno.\else \fi\the\thmno.~{\bi Notation.} }\,}
\def\qst{\advance\thmno by1
\bigpagebreak\ni{\bf
\if\setsections Y\the\scno.\else \fi\the\thmno.~{\bi Question.} }\,}

\def\apf{approximate positive factorization}

\sect{Introduction}

In recent years there has been some interest in the characterization of
the \op s on a \hs\ that admit factorizations
into products of ``nice'' operators, such as normal, selfadjoint, or
positive \op s. Quite different results are obtained depending on
whether
the dimension of the \hs\ is finite or infinite. In the \fd\ case
every \op\ is a product of two normal \op s, by polar
decomposition; an \op\ is a product of selfadjoint \op s \ifo\ its
determinant is real (H.~Radjavi~[\RAD]); and an \op\ is a product
of positive \op s \ifo\ its
determinant is nonnegative (C.S.~Ballantine~[\BAL]).
The infinite-dimensional case was analyzed by
P.Y.~Wu in~[\WU]. There it is shown that for \op s on
a separable in\fd\ \hs, the products of normal \op s, of
selfadjoint \op s, and of positive \op s all form the same class.
Moreover, an \op\ belongs to this class \ifo\ it is a norm limit of
invertible \op s.

The study of factorization into positive and selfadjoint operators
has been extended to several classes of \ca s: homogeneous \ca s
(N.C. Phillips [\PhC]), and unitized stable \ca s and purely
infinite simple \ca s (M. Leen [\Leen]).

The question of approximate positive factorization
was first raised for operators on Hilbert space, by
M. Khalkali, C. Laurie, B. Mathes, and H. Radjavi in [\KLMR].
This question has has been studied in a \ca\ context by T.~Quinn.
In~[\QUI] he shows that if $A$ is an
AF-\al, then every \el\ of $A$ is a norm limit of products of
positive \el s of
$A$ \ifo\ $A$ admits no nonzero \fd\ quotient \ca.

In this paper we pursue the question of approximate positive
factorization further.
Our approach is  new, involving additive commutators. Moreover, our
emphasis is on a global
as opposed to a local analysis; that is, we consider the question
of determining global conditions on an \al\ that ensure that
{\it all} \el s are approximable by products of positive \el s.
In this case, the \al\ is said to have the \apf\ property. (This is the
definition for the unital case; that for the nonunital case is
different. See below.) Our results
lead to new classes of \ca s having the \apf\ property---for example,
type~II$_1$ factors, and
infinite-dimensional simple unital direct limits with slow dimension
growth, real rank zero, and trivial $K_1$~group.

C*-algebras having the \apf\ property are shown to have some nice
properties. For instance, it is shown in the unital case
that the $K_0$~group separates the tracial states.

Our analysis leads to the introduction of a new
concept of rank for a \ca\ that may be of interest in the future.\qedd

The paper is organized as follows: In \S2 we derive some properties of
\ca s having the \apf\ property. In~\S3 we consider
some operations under which this class of \ca s is closed. In the final
section, \S4, we construct large classes of examples of \ca s with the
\apf\ property.\qedd

\sect{The \APFP}

We begin by setting up some notation.

\ntn
If $A$ is a unital \ca , then
$\Inv (A)$ is the invertible group of $A$, and
$\Invo (A)$ is its identity component. Also, $U(A)$ is the unitary group
of $A$ and $U_0 (A)$ is its identity component.

If $A$ is a \ca , then $\tilde{A}$ denotes
its unitization (with a new unit adjoined even if $A$ already has one),
and
$A + 1$ denotes $\{a \in \tilde{A} : a - 1 \in A\}$.
If $A$ is nonunital, we define
$\Inv (A) = \Inv (\tilde{A}) \cap (A + 1)$ and
$\Invo (A) = \Invo (\tilde{A}) \cap (A + 1)$, and
define $U (A)$ and $U_0 (A)$ similarly. \qedd

If $A$ is unital, then there is a canonical isomorphism
$\tilde{A} \cong A \oplus \cn.$ If in this case we identify an
invertible \el\  $a \in A$ with $(a, 1) \in A \oplus \cn,$ then
the definitions given for the nonunital case,
when applied to unital \ca s, agree with the ones
for the unital case.

\dfn
Let $A$ be a unital \ca . Let $P(A)$ denote
the set of finite products of positive \el s of $A$. We say that $A$ has
the {\sl \APFP}  if $P(A)$ is dense in $A$. If $A$ is not unital, then
we say $A$ has the {\sl \APFP}  if the set $P(A + 1)$
of finite products of positive elements in $A + 1$ is dense in $A + 1$.
\qedd

If $A$ is unital, one easily checks that $P(A)^- = A$
\ifo\  $P(A+1)^- = A + 1.$
Thus, we can (and will)
prove results for both the unital and nonunital cases using only the
setup for the nonunital case.

It will be clear from results below that a unitization
$\tilde{A}$ can never have the \APFP .
It is, however,
not at all clear what the relation is between our definition for
nonunital \ca s and the
condition that products of positive elements of $A$ be dense in $A$.
We use the definition above because it works in our theorems.

We make a few comments on $P(A)$.
Obviously, it is multiplicatively closed. Also,
it is invariant under similarity; that is, if $v$ is an invertible
\el\ of $A$, then $vP(A)\inv v= P(A).$ (This is easily
seen by polar decomposing $v$ into a product of a unitary \el\ and a
positive \el.)  Finally, we note that one sided invertibility implies
invertibility for positive elements, so that if
$a = a_1 \cdots a_n \in P(A)$, with $a_1, \dots, a_n$ positive,
and if $a$ is invertible, then so are the $a_j$.

In Theorem 2.4 below, we give some consequences of the \APFP .
We need the following lemma. It is the analog of the fact that if $A$
has stable rank $1,$ then so do all matrix algebras over $A.$
We denote by $M_n (A)$ the algebra of $n \times n$ matrices over
$A,$ and write $M_n$ when $A = \cn.$

\lem
Let $A$ be a unital \ca\ in which $\Invo (A)$ is dense. Then
$\Invo (M_n(A))$ is dense in $M_n(A)$ for each positive integer~$n$.

\pf
We prove the result by induction on~$n$. If $x$ is an \el\ of
$M_n(A)$, we may write it in the form
$$x=\pmatrix a&b\\
             c&d\endpmatrix,$$
where $a\in M_{n-1}(A)$, $d\in A$ and $b$ and $c$ are appropriate column
and row matrices with entries in~$A$.
To show that $x\in (\Invo (M_n(A)))^-$, we may suppose that $d$
belongs to $\Invo (A)$, since,
by assumption, $\Invo (A)$ is dense in $A$.
Hence, we may set $y=a-b\inv dc$. Then
$$x=\pmatrix 1&b\inv d\\ 0&1\endpmatrix\pmatrix y&0\\ 0&d\endpmatrix
\pmatrix 1&0\\ \inv dc&1\endpmatrix.$$
The second factor is in  $(\Invo (M_n(A)))^-$
by the inductive
hypothesis, since $y\in (\Invo (M_{n-1}(A)))^-$.
The other two factors are
clearly in $\Invo (M_n(A))$. So $x \in (\Invo (M_n(A)))^-$. \qed

\thm
Let $A$ be any \ca\  (unital or not) with the \APFP . Then:

(1) $\Inv (M_n (A))$ is connected for all $n$.

(2) $K_1 (A) = 0$.

(3) The (topological) stable rank $\oper{sr} (A)$ is $1$.

(4) The connected stable rank $\oper{csr} (A)$  ([\RfA], Definition 4.7)
is $1$.

\pf
To avoid repetition, we use the setup for the nonunital case.

Functional calculus shows that every positive element in $A + 1$ is a
limit of positive invertible elements in $A + 1$. Such elements are
clearly in $\Invo (A)$. Therefore $\Invo (A)$ is dense in $A + 1$.
By factoring out scalars, we see that $\Invo (\tilde{A})$
is dense in $\tilde{A}$.
This proves (3), and (1) now follows from the previous lemma.  (Note
that $\Inv (B)$ is connected \ifo\  $\Inv (\tilde{B})$ is connected.)
Conclusion (2) is immediate from (1). For (4), we note that, by the
remarks after Corollary 4.10 of~[\RfA], it suffices to show that
left invertible elements in $A$ are invertible and that $\Inv(A)$
is connected. The first statement follows from (3) and the second
from (1).
\qed

If $A$ is any \uca, we denote by $G(A)$ the set of all \el s
$\l\in \tn$, where $\tn$ denotes the set of unit-modulus scalars,
\st\ $\l 1\in P(A)^-$. Clearly, $G(A)$ is a closed subgroup of $\tn$,
and
\tf\ it is the group of all $n$th roots of unity, for some positive
integer~$n$, or it is equal to~$\tn$.

A number of concepts of rank for \ca s have been introduced in recent
years, such as the real, analytic, and stable ranks.
(See [\BP,\MUR,\RIE].)
Using the group $G(A)$ we can define another
rank function. It will be technically useful here, and may be of
independent interest in the future.

\dfn
We set $\rk(A)=n$ if $G(A)$
is finite with $n$~\el s and we set $\rk(A)=\i$ if $G(A)=\tn$.
\qedd

\pro
The function $\rk$ has the following properties for unital \ca s.

(1) $\rk (M_n) = n$.

(2) The algebra $B(H)$ of bounded operators on a separable
infinite-dimensional Hilbert space $H$ satisfies $\rk (B(H)) = \infty$.

(3) If $\f : A \to B$ is a unital \hm , then $\rk (A)$ divides
$\rk (B)$. (Any integer is considered to divide $\infty$.)
In particular, if $B$ is a unital \ca\  and $A$ is a \cs\  of
$B$ containing the identity of $B$,  then $\rk (A)$ divides
$\rk (B)$.

(4) If $A$ has a nonzero representation on $\cn^n$, then
$\rk (A) \leq n$.

(5) If $\rk (A) > 1$, then the commutator ideal of $A$ is all of $A$.

\pf
Parts (1) and (2) follow from Ballantine's theorem~[\Bll] and Wu's
characterization of $P(B(H))$ [\WU] respectively. For (3), note that
$\f (G(A))$ is a subgroup of $G(B)$. Part (4) is immediate from part
(3), and part (5) follows since if the commutator ideal is proper,
then $A$ has a one-dimensional representation. \qed

If $X$ is a compact Hausdorff space, then $M_n \otimes C(X)$ has
rank $n.$ This follows from [\PhC], and is also easy to check
directly.

We also observe that a unital \ca\ $A$  with infinite rank can never
be of type I. Indeed, if $A$ is unital and of type I, choose a maximal
ideal $M$ of $A$, and note that $A/M \cong M_n$ for some $n < \infty.$
Thus, $\rk (A) \leq n.$

A unital \ca\  with the \APFP\  has infinite rank, so we obtain:

\pro
Let $A$ be a unital \ca\  with the \APFP . Then:

(1) $A$ has no nonzero finite-dimensional representations.

(2) The commutator ideal of $A$ is all of $A$.
\rm

\medpagebreak

We next establish a connection between the \APFP\  and the behavior
on the $K_0$-group
of tracial states. We will need part (2) of the following lemma;
part (1) will be used at the end of this section.

\lem
Let $A$ be a unital \ca , let $\t : A \to \cn$ be a selfadjoint
(but not necessarily positive) bounded linear tracelike functional, and
let
$a \in A$ be selfadjoint. Then:

(1) If $\exp (2 \pi i a) \in P(A)$ then $\t (a) \in\t_* (K_0 (A)).$

(2) If $\exp (2 \pi i a) \in P(A)^-$ then $\t (a) \in \t_* (K_0 (A))^-.$

\pf
The proofs of the two parts are similar, but the proof of (2)
is slightly more complicated. Therefore we only prove (2).

Clearly we may assume $\t \neq 0.$ Let $\e > 0.$
Use polar decomposition
to choose positive elements $b_1, \dots, b_n \in A$
such that $z = b_1 \cdots b_n$ is unitary and
$\n{z - \exp(2 \pi i a)} < 2 \pi \e / \n{\t}.$ Set
$h = \frac{1}{2 \pi i} \log (\exp (2 \pi i a)^* z)$. Then $h$ is
selfadjoint and satisfies
$z = \exp(2 \pi i a) \exp(2 \pi i h)$ and $\n{h} < \e / \n{\t}.$
Let $\dd$ be the de la Harpe-Skandalis
determinant~[\dlHS] associated to $\t;$ it is a \hm\  from
$\Invo(A)$ to the additive group $\cn / \t_* (K_0 (A)).$
Proposition 2(b) of~[\dlHS] states that
$\dd (\exp(c)) = \frac{1}{2 \pi i} \t (c) + \t_* (K_0 (A))$
for any $c \in A$. Computing $\dd (z)$ in two different ways,
we obtain
$$ \frac{1}{2 \pi i} \sum_{k = 1}^n  \t (\log(b_k)) - \t (a) - \t (h)
                     \in \t_* (K_0 (A)).$$
(Note that $\log(b_k)$ is defined since $b_k$ is invertible.)
Thus, we can write
$$ \frac{1}{2 \pi i} \sum_{k = 1}^n  \t (\log(b_k)) =
     \t (a) + \t (h) - \alpha $$
for some $\alpha \in \t_* (K_0 (A)).$
The right hand side of
this equation is real and the left hand side is purely imaginary.
Therefore both sides are zero.
Since $\m{ \t (h)} \leq \n{\t} \n{h} < \e,$ we get
$\m{ \t (a) - \alpha} < \e.$ Thus
$\oper{dist} (\t (a), \t_* (K_0 (A))) < \e.$
Since $\e > 0$ is arbitrary, this shows that
$\t (a)  \in \t_* (K_0 (A))^-.$ \qed

As a corollary, we obtain the following theorem.

\thm
Let $A$ be a unital \ca\  with the \APFP . Then $K_0 (A)$ distinguishes
the tracial states on $A$. That is, if $\t_1$ and $\t_2$ are two
tracial states on $A$ such that $(\t_1)_* = (\t_2)_*$ as maps from
$K_0 (A)$ to $\rn$, then $\t_1 = \t_2$.

\pf
Suppose $K_0 (A)$ does not distinguish the tracial states on $A$.
Let $\t_1$
and $\t_2$ be two distinct tracial states such that
$(\t_1)_* = (\t_2)_*$. Set $\t = \t_1 - \t_2.$ Then
$\t_* (K_0 (A)) = \{0\}.$ Since $\t_1 \neq \t_2,$ there is
a selfadjoint $a \in A$ such that $\t_1 (a) \neq \t_2 (a).$
Then $\exp (2 \pi i a) \not\in \CPP{A}$ by part (2) of the previous
lemma.  \qed

For some classes of \ca s (see~[\BDR] and Theorem 4.9 below), the
conclusion of this theorem implies that the algebra has real rank
zero~[\BP]. Example 4.16 below will show, however, that the
\APFP\  does not in general imply real rank zero.

Lemma  2.8(1) yields results on $P(A)$ which show why one should
consider the condition $\CPP{A} = A$ rather than $P(A) = A.$

\pro
Let $A$ be a separable unital
\ca\  which has a tracial state $\t.$
(Equivalently, the linear span of the commutators $ab - ba$ is not
dense.) Then there is $\l \in \cn$ such
that $\l \cdot 1 \not\in P(A).$

\pf
Since $A$ is separable, $K_0 (A)$ is countable, so
$\t_* (K_0 (A))$ is a proper subset of $\rn.$ Choose
$\alpha \in \rn$ such that  $\alpha \not\in \t_* (K_0 (A)).$
Then $\exp ( 2 \pi i \alpha ) \cdot 1 \not\in P(A)$ by Lemma  2.8(1).
\qed

\rem
This proposition shows that if $A$ is a simple
separable unital stably finite
exact \ca , then there is $\l \in \cn$ such
that $\l \cdot 1 \not\in P(A).$ Indeed,
it follows from Theorem 6.1 of~[\RrB] that a simple stably finite
\ca\  has a nontrivial
quasitrace. By [\Hg],  every quasitrace on an exact \ca\  is a trace.
\qedd

We can even use Lemma  2.8 to shed some more light on the rank
defined above:

\pro
Let $A$ be a unital \ca . Suppose $n$ is a positive integer and $\t$ is
a tracial state on $A$ such that $n \t_* (K_0 (A)) \subset \zn.$
Then $\rk (A) | n.$

\pf
The hypotheses imply that $\t_* (K_0 (A))$ is closed in $\rn.$
Lemma  2.8(2) therefore implies that if $\alpha \in \rn$ and
$\exp ( 2 \pi i \alpha ) \cdot 1 \in P(A)^-,$ then $n \alpha \in \zn.$
Thus $G(A)$ is a subgroup of $\{ \zt \in \tn : \zt^n = 1\}.$ \qed

\pro
If $A$ is unital and has infinite rank, then $\t_* (K_0 (A))$ is
dense in $\rn$ for any tracial state $\t.$

\pf
This is immediate from Lemma  2.8(2). \qed

This proposition shows that if $A$ has a tracial state and has infinite
rank,
then $K_0 (A)$ must be large. Note, however, that the Cuntz algebra
${\cal O}_2$ has infinite rank (since, by [\Leen],
$\Invo (A) \subset P(A)$ for purely infinite simple unital \ca s $A$),
but $K_0 ({\cal O}_2) = 0.$

\exa
Let $C^*_{\oper{r}} (F_m)$ denote the reduced \ca\  of the free group
on $m$ generators. Then $\rk (M_n \otimes C^*_{\oper{r}} (F_m)) = n.$
Indeed, we obtain $\rk (M_n \otimes C^*_{\oper{r}} (F_m)) \geq n$
by Proposition 2.6(1) and (3), applied to the subalgebra
$M_n \otimes 1.$ The reverse inequality follows from
Proposition~2.12.

The same argument shows that if $A$ is any unital
\ca\  with a tracial state $\t$ for which $\t_* (K_0 (A)) = \zn,$
then $\rk (M_n \otimes A) = n.$
\qedd

This example gives simple infinite-dimensional \ca s $A$ with
arbitrary finite values of $\rk (A).$ On the other hand,
Theorem~4.10 below gives many stably finite simple unital \ca s
with infinite rank.

\sect{Invariance properties}

In this section, we prove that the \APFP\  is preserved under several
natural operations on \ca s. Our first goal is to show that if
$A$ has the \APFP , then so does $M_n (A).$ The following notation
will be useful.
To avoid repetition,  we work in the setup for the
nonunital case.

\ntn
Let $A$ be a \ca.
We denote by $[A,A]$ the closed linear span of
the commutators $[a,b] = ab - ba$ of $A$ and by $L(A)$ the set of all
\el s $a$ of
$A$ \st\ $\exp ({\l a})$ belongs to $P(A + 1)^-$ for all $\l\in \cn$.
\qedd

Note that $[\tilde{A}, \tilde{A}] = [A, A]$, and
that if $a \in A$ then $\exp (a)$ is always in $A + 1$.
Further, if $A$ is unital, $a \in A,$ and
$\exp (a),$ calculated in $\tilde{A},$ is in $P(A+1)^-,$
then $\exp (a),$ calculated this time in $A$, is in $P(A)^-.$
Also, it is obvious that if
$L(A)=A$ and $\Invo (\tilde{A})$ is dense in $\tilde{A}$,
then $P(A + 1)$ is dense in $A + 1,$ and conversely.

Before proceeding, we need to recall some known
results about exponentials.
If $a,b\in A$, then
$$\exp ({a+b}) =\limn (\exp ({a/n}) \exp ({b/n}) )^n \tag 1$$
and
$$\exp ({[a,b]}) = \limn(\exp ({-a/n})
            \exp ({-b/n}) \exp ({a/n}) \exp ({b/n}) )^{n^2}.\tag 2$$
For lack of a better source, we
refer to the proof of Lemma 2.4 of [\PhB].
\qedd

\lem
Let $A$ be a \ca. Then $L(A)$ is a closed, linear subspace of
$A$ containing $[A,A]$.

\pf It is obvious that $L(A)$ contains $0,$ is norm closed, and is
closed under
multiplication by scalars. That $L(A)$ is closed under addition follows
from Equation~(1). To see that $[A,A]\ss L(A),$ it suffices to show that
$[a,b]\in L(A)$ for all $a,b\in A$. Hence, it is sufficient to show that
$\exp ({[a,b]}) \in P(A + 1)^-$ if $a, b\in A$ and $b=b^*$. For such
\el s
$a,b$ we have, by Equation~(2),
$$\exp ({[a,b]}) =\limn(v_n \exp ({-b/n})\inv v_n \exp ({b/n}))^{n^2},$$
where $v_n=\exp ({-a/n})$. Since $\exp ({b/n})$ is positive and since
$P(A + 1)$ is invariant under
similarity, we have $v_n \exp ({-b/n})\inv v_n \exp ({b/n})\in P(A + 1)$
and
\tf\ $\exp ({[a,b]})\in P(A + 1)^-$, as required.\qed

\exa  Let $H$ be an infinite-dimensional \hs, and let $A = B(H).$
Since every \op\ on $H$ is
a sum of two commutators ([\HAL], Corollary 2 to Problem 234),
we have $[A,A]=A$.
Hence, $\exp (a) \in P(A)^-$
for all $a\in A$, by the previous lemma. Of course,
$\Invo (A) =\Inv (A)$ in this case, so $P(A)^-
=(\Inv (A))^-$. This generalizes, for if $A$ is any properly infinite
von~Neumann \al, then $[A,A]=A$ [\PT] and \tf, $P(A)^-=(\Inv (A))^-$,
again by the previous lemma.\qedd

\thm
Let $A$ be a \ca\ having the \apf\ property. Then for each positive
integer~$n$, the matrix \al\ $M_n(A)$ also has the \apf\ property.

\pf
We first show that if $a=(a_{ij})$ is an \el\ of $M_n(A)$ with
zero diagonal, then $a$ is a limit of sums of
commutators. To do this we identify $M_n(A)$ with the tensor product
$M_n \otimes A$ in the usual way, so that, if the \el s
$e^{ij}$ ($i,j=1,\dots,n$)
form the canonical basis for $M_n$, we may write
$a=\sum_{i,j} {e^{ij}\otimes a_{ij}}$. Let $(e_{\l})$ be an
approximate identity for $A$. Then the equation
$$a = \lim_{\l} \sum_{i, j}
               [e^{ij}\otimes a_{ij},e^{jj}\otimes e_{\l}]$$
expresses $a$ as a limit of sums of commutators.

Now let $a$ be an arbitrary \el\ of $M_n(A)$. Clearly,
$a=b+c$, where $b$ has zero diagonal and $c$ is a diagonal matrix.
Since $b$ is a limit of sums of commutators, $b$ belongs to $L(M_n(A))$
by Lemma~3.2, and since $L(A)=A$ it is
clear that $c$ belongs to $L(M_n(A))$. Hence, $L(M_n(A))=M_n(A)$.
Moreover, since
$(\Invo (\tilde{A}))^- = \tilde{A}$, we have
$(\Invo (M_n(\tilde{A})))^- = M_n(\tilde{A})$ by Lemma~2.3.
Combining these two facts we get $P(M_n(A) + 1)^- = M_n(A) + 1$.\qed

\thm
(1) The finite direct sum of \ca s with the \APFP\  has the \APFP .

(2) The direct limit of \ca s with the \APFP\  has the \APFP .

(3) If $A$ has the \APFP\  and $F$ is a finite-dimensional \ca, then
$A \otimes F$  has the \APFP .

(4) If $A$ has the \APFP\  and $F$ is an AF algebra, then
$A \otimes F$  has the \APFP .

\pf
Parts (1) and (2) are immediate. Part (3) follows from (1) and the
previous theorem, and part (4) follows from (2) and (3). \qed

\thm
Let $A$ be a \ca , and let $I$ be a closed ideal in $A$. Then:

(1) If $I$ and $A/I$ have the \APFP , then so does $A$.

(2) If $A$ has the \APFP , then so does $A/I$.

\pf
Let $B = A/I$ and let $\pi : \tilde{A} \to \tilde{B}$ be the quotient
map.

(1) Theorem 2.4(3) and (4),  and Theorem 4.11 of~[\RfA],  imply that
$\oper{sr} (A) = 1$. By polar decomposition, it therefore suffices to
prove that $U(\tilde{A}) \cap (A + 1) \subset \CPP{A + 1}$.
Let $u \in U(\tilde{A}) \cap (A + 1)$, and
let $\e > 0$. Continuity of the polar decomposition allows us
to use the \APFP\  in $B$ to find a unitary $y \in P(B + 1)$
such that $\n{y - \pi (u)} < \min (1, \e / 2)$.

Write $y = b_1 \cdots b_m$ with $b_1, \dots, b_m \in B + 1$ positive.
Then the $b_j$ are also invertible. Hence for each $j$ there
is $\r_j > -1$ such that $\oper{sp} (b_j - 1) \subset [\r_j, \infty)$.
Using functional calculus, it is easy to choose $a_j \in A$ selfadjoint
such that $\pi (a_j) = b_j -1$ and
$\oper{sp} (a_j) \subset [\r_j, \infty)$. In particular, $a_j + 1$ is
invertible. Set
$x_0 = (a_1 + 1) \cdots (a_m + 1)$, and let $x$ be the unitary part of
the polar decomposition of $x_0$. Then $x \in P(A + 1)$ and
$\pi (x) = y$. Therefore $\n{\pi (x^* u - 1)} < \e / 2$. Clearly also
$x^* u - 1 \in A$.

Choose $c_0 \in A$ such that $\pi (c_0) = \pi (x^* u - 1)$ and
$\n{c_0} < \e / 2$. Set $c = x^* u - c_0$. Then $\pi (c) = 1$,
so that $c \in I + 1$. Therefore there is $r \in P(I + 1)$ such that
$\n{r - c} < \e / 2$. It follows that $xr \in P(A + 1)$ and
$\n{xr - u} \leq \n{r - c} + \n{c_0} < \e$.

(2) This is immediate from the fact that $\pi (A + 1) = B + 1$.
\qed

We will see below (Example 4.16) that the \APFP\  does not always
pass to ideals.

\sect{Characterizations and classes of examples}

In this section, we construct a variety of examples of \ca s having the
\APFP . In some special classes we can even give necessary and
sufficient conditions for the \APFP .

We begin by recalling the following theorem of Quinn~[\QUI]:

\thm
A unital AF algebra has the \APFP\  \ifo\  it has no nonzero
finite-dimensional quotient \ca s.
\rm

\medpagebreak

Next, we consider von Neumann factors. We need two lemmas.

\lem
Let $A$ be a unital \ca .
Then $[A,A]+\cn 1=A$ \ifo\  $A$ has at most one tracial state.

\pf
This follows from Proposition 2.7 of~[\CP], using Theorem 2.6 and the
definitions in Section 2 of the same paper.
\qed

A simpler proof than that of~[\CP] can be obtained as follows. The
forward implication is obvious. As for the reverse implication,
if $[A, A] + \cn 1 \neq A,$ then the
Hahn--Banach theorem can be used to find a nonzero selfadjoint bounded
linear functional
which vanishes on $[A,A]+\cn 1$, and the Jordan decomposition
then yields two different tracial states. \qedd

\lem
Let $A$ be a unital \ca\ \st\ $[A,A]+\cn 1=A$ and \st\ $\Invo (A)$ is
dense in~$A$. Then for each \el\ $a\in A$ \tes\ $\l\in \tn$
\st\ $\l a\in P(A)^-$.

\pf If $a\in A$, then for some scalar $\l$ the \el\ $a-\l$ belongs to
$[A,A]$ and \tf\ $\exp ({a-\l})$ belongs to $P(A)^-$, by Lemma~3.2.
Hence,
$\l'\exp (a) \in P(A)^-$, where $\l'=\exp ({-\l}) /\m{\exp ({-\l}) }$
and $\m{\l'}=1$.
It is an immediate consequence that for each \el\ $a$ of
$\Invo (A)$ \tes\ $\l\in \tn$ \st\ $\l a\in P(A)^-$.
That the same result \tf\ holds for any \el\ $a$ of
$A$ now follows from the assumption that $\Invo (A)$ is dense in $A$.
\qed

\thm
A factor has the \apf\ property \ifo\  it is of type~II$_1$.

\pf Let $A$ be a factor with the \APFP . Then $sr(A) = 1$, so $A$
is finite. Proposition~2.7(1) implies that $A$ is infinite-dimensional.
Hence, $A$ must be of type~II$_1$.

Now suppose conversely that $A$ is of type~II$_1$. Then
it is simple and has a unique tracial state, since, as is well
known, this is true of any finite factor. Consequently, by
Lemma~4.2, $[A,A]+\cn 1=A$. Moreover, by a result of
H.~Choda~[\CHO], finiteness of $A$ implies $\Invo (A)$ is dense in $A$.
Hence, by the previous lemma, for each \el\ $a$ of $A$,
\tes\ a scalar~$\l$ of unit modulus \st\ $\l a$ belongs to $P(A)^-$.
To complete the proof it is clear that we need now show only that
$\rk(A)=\i$.  It is well known that a type~II$_1$ factor contains a
unital copy of $M_n$ for any $n$, and it therefore
follows from Ballantine's theorem [\BAL]
that $G(A)$ contains all $n$th roots of unity for every $n$.
Therefore $G(A) = \tn$.
\qed

\exa
The previous theorem, combined with the method of Example 3.10 of~[\Ph],
yields a simple separable nonnuclear (hence not AF) unital
\ca\  with the \APFP .
\qedd

\thm
Let $F$ be an in\fd\ UHF \al\ and let $A$ be a unital \ca\  with
at most one tracial state. Suppose either

(1) $\Invo (A)$ is dense in $A$, or

(2) $A$ is simple, $K_1 (A) = 0,$ and $A$ has a tracial state.

\noindent
Then the C*-tensor product $F\otimes A$ has the \apf\ property.

\pf
Let $C = F \otimes A.$
Write $F = \dirlim M_{k(n)}$ (with unital maps), and set
$C_n = M_{k(n)} \otimes A.$ Thus $F \otimes A = \dirlim C_n.$
In case (1), $\Invo (C_n)$ is dense in $C_n$ by Lemma 2.3,
and it follows
that $\Invo (C)$ is dense in $C.$ In case (2), $C$ has
stable rank $1$ by Corollary 6.6 of~[\RrA]. Also, $K_1 (C) = 0,$ so
$\Invo (C) = \Inv (C)$ by Theorem~2.10 of~[\Rf]. Therefore $\Invo (C)$
is dense in $C.$
Thus, in either case, $\Invo (C)$ is dense in $C.$

In both cases, each $C_n$ again has at most one tracial state, so
the same holds for $C$, and Lemma~4.2 implies that
$[C, C] + \cn 1 = C.$ Lemma~4.3 now implies that for all $c \in C$
there is $\l \in \tn$ such that $\l c \in \CPP{C}.$
To prove the theorem, it remains only to prove that
$\inv{\l}$ is in $\CPP{C}.$ But
it follows from Proposition 2.6(1) and (3) that $\rk (F) = \infty,$
whence also $\rk (C) = \infty.$ Thus, $\tn \cdot 1 \subset \CPP{C},$
and the theorem follows. \qed

\exa
In part (1) of the previous theorem, we can take $A$ to be $\cn$
(giving a fast proof that infinite-dimensional
UHF algebras have the \APFP , a special case
of Quinn's results as stated in Theorem 4.1 above), a
type~II$_1$ factor,
or the unitization of any stable \ca\  $B$ satisfying $K_1 (B) = 0$
and $\oper{sr} (B) = 1.$
In part (2), we can take $A$ to be the Choi algebra
$C^*_{\oper{r}} (\zn / 2 \zn * \zn / 3 \zn).$
(See Example (2) in Section 3 of~[\Cu] and Theorem 2.8 (and its proof)
in~[\Choi] for verification
that the hypotheses of (2) are satisfied.)

We note that Theorem 7.2 of~[\RrB] shows that the algebras
$F \otimes A$ in (2) have real rank zero whenever every quasitrace is
a trace. This happens for all exact \ca s [\Hg], and thus in particular
for all nuclear \ca s.
\qedd

We are now going to give necessary and sufficient conditions
for an infinite-dimensional
simple direct limit $A$ with slow dimension growth as in~[\BDR] to
have the \APFP . Recall that slow dimension growth
means that $A$ is a unital direct
limit $\dirlim A_i,$ in which each
$A_i$ is of the form
$A_i = \bigoplus_{t = 1}^{s(i)} C(X_{it}) \otimes M_{n(i, t)}$
for finite-dimensional connected compact
Hausdorff spaces $X_{it},$ such that
$$
\lim_{i \to \infty} \sup_{1 \leq t \leq s(i)}
                      \oper{dim} (X_{it}) / n(i, t) = 0.
$$
It is implicitly assumed (but not stated) in~[\BDR] that the algebras
are supposed to be infinite-dimensional. Note (see~[\BDR]) that if $A$
is as above, and is simple, then
$$
\lim_{i \to \infty} \inf_{1 \leq t \leq s(i)} n(i, t) = \infty.
$$

The following lemma will be needed.

\lem
For every $\e > 0$ there is $N \in \nn$ such that whenever
$n \geq N$ and $\l_1, \dots, \l_n \in \tn$, then
there exist $\mu_1, \dots, \mu_n \in \tn$ with $\mu_1 \cdots \mu_n = 1$
and such that the set $\{\mu_1^{-1} \l_1, \dots, \mu_n^{-1} \l_n\}$
is $\e$-dense
in $\tn$. (That is, given $\l \in \tn$, there is $k$ such that
$| \l - \mu_k^{-1} \l_k | < \e$.)

\pf
Choose $m$ such that the primitive $m$-th root of $1$ given by
$\zt = \exp( 2 \pi i / m)$ satisfies $| \zt - 1| < \e / 2$.
Let $N = m^2$. Let $n \geq N$ and let  $\l_1, \dots, \l_n \in \tn$.
The points $1, \zt, \dots, \zt^{m - 1}$ divide $\tn$ into $m$
arcs, and it follows that one of these arcs, say the one from
$\zt^k$ to $\zt^{k + 1}$, must contain at least $m$ of the points
$\l_j$. Renumbering, we may assume that $\l_1, \dots, \l_m$ are in this
arc. Set $\mu_1 = 1,\ \mu_2 = \zt, \  \dots, \mu_m = \zt^{m - 1}$, and
let $\mu_{m + 1}, \dots, \mu_n$ be arbitrary elements of $\tn$
satisfying $\mu_{m + 1} \cdots \mu_n = \zt^{-m(m - 1)/2}$. Then
$|\mu_j^{-1} \l_j - \zt^{k - j + 1}| < \e / 2$ for $1 \leq j \leq m$,
and $\{\zt^k, \zt^{k - 1}, \dots, \zt^{k - m + 1} \}$ is $\e / 2$-dense
in $\tn$, so  $\{\mu_1^{-1} \l_1, \dots, \mu_n^{-1} \l_n\}$ is
$\e$-dense.
\qed

\thm
Let $A$ be an infinite-dimensional simple unital direct limit with slow
dimension growth as
in~[\BDR]. Then $A$ has the \APFP\  if and only if $A$ has real rank
zero and $K_1 (A) = 0$.

\pf
The condition $K_1 (A) = 0$ is necessary by Theorem~2.4(2).
To see that real rank zero is necessary, first note
that $K_0 (A)$ distinguishes the tracial states by Theorem 2.9.
The Riesz decomposition property (Theorem 2.7 of~[\Gd]) then shows
that the projections in $A$ distinguish the tracial states.
Theorem 2 of~[\BDR]
now implies that $A$ has real rank zero.

To prove sufficiency, note that Theorem~1 of [\BDR] implies that
$\oper{sr} (A) = 1$. Theorem 2.10 of~[\Rf] therefore implies
that $\Inv (A) / \Invo (A) \to K_1 (A)$ is an isomorphism. Thus,
the hypotheses imply that $\Invo (A)$ is dense in $A$.

By polar decomposition, it now suffices to show that
$U_0 (A) \subset  \CPP{A}$. Write $A = \dirlim A_i$, where each
$A_i$ is of the form
$A_i = \bigoplus_{t = 1}^{s(i)} C(X_{it}) \otimes M_{n(i, t)}$,
in such a way that the system has slow dimension growth and all the
maps are unital. Let $\f_i : A_i \to A$ and
$\f_{i_1, i_2} : A_{i_1} \to A_{i_2}$ be the associated maps.
Without loss of generality we may assume that the
map from each $C(X_{it}) \otimes M_{n(i, t)}$ to $A$ is nonzero.
Let $u \in U_0 (A)$ and let $\e > 0$.

The first thing we have to do is ensure that we are working with
elements with large spectrums. Choose $n_0$ to be the $N$ that works in
the previous lemma for $\e / 8$ in place of $\e$. Choose $i_0$ so large
that some $n(i_0, t_0)$ is at least $n_0$ and there is
$u_0 \in U_0 (A_{i_0})$ with $\n{u - \f_{i_0} (u_0)} < \e / 8$.
The kernel of the
map from $C(X_{i_0, t_0}) \otimes M_{n(i_0, t_0)}$ to $A$ has the form
$C_0 (V) \otimes M_{n(i_0, t_0)}$ for some proper open subset
$V \subset X_{i_0, t_0}$.  Let $x_0 \in X_{i_0, t_0}$ be any point not
in $V$. Then $u_0 (x_0)$ makes sense as
an element of $M_{n(i_0, t_0)}$, and in a suitable basis can be
written as $u_0 (x_0) = \oper{diag} (\l_1, \dots, \l_{n(i_0, t_0)})$.
Choose $\mu_1, \dots, \mu_{n(i_0, t_0)}$ according to the conclusion of
the previous lemma. Define $v_0 \in U (A_{i_0})$ to be the constant
function $1$ on $X_{i_0, t}$ for $t \neq t_0$, and the constant
function $\oper{diag} (\mu_1, \dots, \mu_{n(i_0, t_0)})$ on
$X_{i_0, t_0}$. Since
$\det (\oper{diag} (\mu_1, \dots, \mu_{n(i_0, t_0)})) = 1$, Ballantine's
theorem~[\Bll] implies that $v_0$ is a product of positive elements. Set
$w_0 = v_0^* u_0$. By considering the value at $x_0 \in X_{i_0, t_0}$,
and noting that $x_0 \not\in V$,
we see that $\oper{sp} (\f_{i_0} (w_0))$ is $\e/8$-dense in $\tn$.

Clearly $\f_{i_0} (w_0) \in U_0 (A)$. Since $A$ has real rank zero,
it follows
from Corollary 6 of~[\Ln] that there is $w_1 \in U_0 (A)$ with finite
spectrum and $\n{w_1 - \f_{i_0} (w_0)} < \e/8$. Thus there are distinct
$\l_1, \dots, \l_m \in \tn$ (not the same as in the previous
paragraph) and nonzero \mops\  $p_1, \dots, p_m \in A$ such that
$w_1 = \sum_{j = 1}^m \l_j p_j$. We take the $\l_j$ to be ordered
cyclically.  The estimates above imply that $\oper{sp} (w_1)$ is
$\e/4$-dense in $\tn$, and it follows that
$\m{\l_j - \l_{j + 1}} < \e/2$ for $1 \leq j \leq m - 1$.

Standard functional calculus arguments yield $i_1$
and \mops\  $q_1, \dots, q_m \in A_{i_1}$ (necessarily
nonzero) such that the unitary $w_2 = \sum_{j = 1}^m \l_j q_j$
satisfies $\n{\f_{i_1} (w_2) - w_1} < \e/8$. Choose $N$ such that
$2 \pi / N < \e / 8$, and use Lemma F of~[\BDR] to choose $i \geq i_1$
such that for every $t$ the image $q_j^{(t)}$ of $\f_{i_1, i}(q_j)$ in
$C(X_{it}) \otimes M_{n(i, t)}$ has rank at least $N + \dim (X_{it})$.
(Note that this rank is constant on $X_{it}$, since $X_{it}$ is
assumed connected.)

For each fixed $t$ we now construct, by induction on $j$,
\mops\  $e_j^{(t)} \in C(X_{it}) \otimes M_{n(i, t)}$
which are trivial in the sense that each is Murray-von Neumann
equivalent to a constant projection, and satisfying
$e_1^{(t)} \leq q_1^{(t)}$, \  $e_j^{(t)}
\leq q_{j - 1}^{(t)} + q_j^{(t)}$ for $j \geq 2$,
$\rk (e_j^{(t)}) \geq N$ for all $j$, and
$\sum_{j = 1}^m e_j^{(t)} = 1$. To make the induction work,
we impose the additional conditions
$$ \sum_{l = 1}^{j} e_l^{(t)} \geq \sum_{l = 1}^{j - 1} q_l^{(t)}
       \,\,\,\,\,\, {\text{and}} \,\,\,\,\,\,
\rk \big(\sum_{l = 1}^{j} q_l^{(t)} - \sum_{l = 1}^{j} e_l^{(t)} \big)
             \leq \dim(X_{it}) / 2. $$

The existence of $e_1^{(t)}$ is immediate from Theorem 2.5(a) of~[\Gd].
Suppose now that $e_1^{(t)}, \dots, e_j^{(t)}$ have been constructed;
we construct $e_{j + 1}^{(t)}$. If $j + 1 < m$, set
$f = \sum_{l = 1}^{j} q_l^{(t)} - \sum_{l = 1}^{j} e_l^{(t)}$.
Theorem 2.5(c) of~[\Gd] implies that $f$ is Murray-von Neumann
equivalent to a subprojection of a trivial projection $g$ of rank
at most $\rk (f) + (\dim(X_{it}) - 1) / 2 \leq \dim (X_{it})$.
We may clearly assume $f \leq g$.
By the same theorem, $g - f$ is Murray-von Neumann equivalent to
a subprojection of $q_{j + 1}^{(t)}$. Thus, we may assume
$f \leq g \leq f + q_{j + 1}^{(t)}$. Part (a) of this theorem now yields
a trivial subprojection $h$ of $f + q_{j + 1}^{(t)} - g$ such that
$$
\rk (f + q_{j + 1}^{(t)} - (g + h)) \leq \dim(X_{it}) / 2.
$$
We take $e_{j + 1}^{(t)} = g + h$. Note that
$$
\align \rk(e_{j + 1}^{(t)}) & \geq
         \rk (f) + \rk(q_{j + 1}^{(t)}) - \dim (X_{it}) / 2  \\
       &  \geq \rk(q_{j + 1}^{(t)}) - \dim (X_{it}) \geq N.
\endalign
$$
This completes the induction step in case $j + 1 < m$. If $j + 1 = m$,
then set $e_{j + 1}^{(t)} = 1 - \sum_{l = 1}^{m - 1} e_l^{(t)}$. Note
that
$$
\rk(e_{j + 1}^{(t)}) \geq N + \dim (X_{it}) \geq \dim (X_{it}) / 2,
$$
and the sum of
$e_{j + 1}^{(t)}$ and the trivial projection
$\sum_{l = 1}^{m - 1} e_l^{(t)}$ is the trivial projection $1$, so that
Theorem 2.5(b) of~[\Gd] implies that $e_{j + 1}^{(t)}$ is trivial.
This completes the induction.

We now define $e_j = \sum_t e_j^{(t)}$ and $w_3 = \sum_j \l_j e_j$.
Since $e_j \leq q_{j - 1} + q_j$, and since
$\m{\l_j - \l_{j - 1}} < \e / 2$, we have
$\n{w_3 -  \f_{i_1, i}(w_2)} < \e / 2$. We will now approximate
$w_3$ by a product of positive elements. For each $j$ and $t$, let
$N_j^{(t)} = \rk (e_j^{(t)}) \geq N$, and let $\a_j^{(t)}$ by the
least nonnegative real number such that
$[\exp (2 \pi i \a_j^{(t)}) \l_j]^{N_j^{(t)}} = 1$. Note that
$0 \leq \a_j^{(t)} \leq 1/N$, so that
$\m{ \exp (2 \pi i \a_j^{(t)}) - 1} \leq 2 \pi / N$. Now define
$w_4 = \sum_{j, t} \exp (2 \pi i \a_j^{(t)}) \l_j e_j^{(t)}$.
Then $w_4$ is a unitary which satisfies
$\n{w_4 - w_3} \leq 2 \pi / N < \e / 8$. Furthermore, for each $j$
and $t$, the choice of $\a_j^{(t)}$ and the fact that $e_j^{(t)}$ is
a trivial projection imply that
$\exp (2 \pi i \a_j^{(t)}) \l_j e_j^{(t)}$ is a constant
scalar multiple of the identity in
$$e_j^{(t)} (C(X_{it}) \otimes M_{n(i, t)}) e_j^{(t)}
    \cong C(X_{it}) \otimes M_{N_j^{(t)}} $$
whose determinant is $1$.  Using Ballantine's theorem on these elements
and forming the direct sum over $j$ and $t$, we therefore see that
$w_4$ is a product of five positive elements.

We now have $\n{u - \f_{i_0} (u_0)} < \e / 8$, \  $u_0 = v_0 w_0$, and
$$\align \n{\f_{i_0} (w_0) - \f_i (w_4)}
       &\leq \n{\f_{i_0} (w_0) - w_1} + \n{w_1 - \f_{i_1} (w_2)}\\
       &\qquad\qquad\qquad+\n{ \f_{i_1, i}(w_2) - w_3} + \n{w_3 - w_4}\\
       &<  \e / 8 + \e / 8 + 4 \e / 8 + e / 8 = 7 \e / 8;
\endalign$$
furthermore, both $v_0$ and $w_4$ are finite products of positive
elements. Thus $u$ can be approximated to within $\e$ by finite products
of positive elements, and the proof that $U_0 (A) \subset \CPP{A}$ is
complete. \qed

This proof actually shows that products of $11$ positive invertible
elements are dense (one from the polar decomposition and five from
each application of Ballantine's theorem).

The proof of this theorem also implies the following result, which holds
regardless of whether the algebra has real rank zero or one.

\thm
Let $A$ be an infinite-dimensional simple unital direct limit with slow
dimension growth as in~[\BDR]. Then the rank of  $A$ is infinite.

\pf
If the element $u$ in the previous proof is a scalar multiple of the
identity, then the unitary $\f_{i_0} (w_0)$ constructed in the proof
will already have finite spectrum. The only use of real rank zero in
that direction of the proof was to replace this unitary with a nearby
one with finite spectrum.
\qed

\exa
Let $X$ be a finite CW complex such that $K^0 (X)$ has
a nontrivial torsion element $\eta$ and $K^1 (X) = 0$.  Choose
a sequence $x_1, x_2, \dots$ in $X$ such that each tail of the sequence
is dense in $X$. Define $A_i = C(X, M_{2^i})$ and define
$\f_{i, i + 1} : A_i \to A_{i + 1} \cong M_2 (A_i)$ by
$\f_{i, i + 1} (f) = \oper{diag} (f, f(x_i) \cdot 1)$.
It follows from~[\GdB] that $A = \dirlim A_i$ is simple and has
real rank zero, and one computes directly that $K_1 (A) = 0$ and
that the image of $\eta \in K_0 (A_0)$ in $K_0 (A)$ is nonzero.
Therefore $A$ is a
separable simple nuclear non-AF algebra, which has the \APFP\  by
Theorem 4.9.

A simple example of such a finite CW complex $X$ can be constructed
as follows. Let $D \subset \cn$ be the closed unit disk. Let $\sim$
be the equivalence relation $z_1 \sim z_2$ if $\m{z_1} = \m{z_2} = 1$
and $z_1^n = z_2^n.$ Then $X = D / \!\! \sim$ is a finite CW complex
such that $K^0 (X) \cong \zn \oplus \zn / n \zn$ and $K^1 (X) = 0.$
\qed

Theorem 4.9 raises the following two questions.

\qst
Let $A$ be a simple infinite-dimensional unital \ca\  with
real rank zero, stable rank one, and $K_1 (A) = 0.$ Does it
follow that $A$ has the \APFP ? \qedd

Theorem 4.9 shows that the answer is yes for infinite-dimensional
direct limits with slow dimension growth.
The answer is also yes if $A$ has infinite rank
and a unique tracial state, by Lemmas 4.2 and 4.3.

\qst
Let $A$ be a simple unital \ca\  with the \APFP . Does it follow
that $A$ has real rank zero? \qedd

Theorem 4.9 shows that the answer to this question is also
yes for infinite-dimensional
direct limits with slow dimension growth.
Furthermore, Theorem 2.9 and Proposition 2.13 imply that if $A$ is a
unital
algebra with the \APFP ,
then at least the stabilization
$K \otimes A$ has a fair number of projections---enough
to distinguish the tracial states,
and, if there is a tracial state $\tau,$ then
enough that $\tau_* (K_0 (A))$ is dense in $\rn.$
The requirements that $A$ be unital and simple are necessary,
as is shown by the next theorem and by Example~4.16.

We now give necessary and sufficient conditions for a stable \ca\  to
have the \APFP .

\thm
Let $A$ be any \ca . Then $K \otimes A$ has the
\APFP\  \ifo\  $\oper{sr} (A) = 1$ and $K_1 (A) = 0$.

\pf
Let $B = K \otimes A$.

If $B$ has the \APFP , then  $\oper{sr} (A) = 1$ by
Theorem~2.4(3)
and Theorem 3.6 of~[\RfA], and $K_1 (A) = 0$ by Theorem 2.4(2).

For the converse, let $\oper{sr} (A) = 1$ and $K_1 (A) = 0$, so that
$\oper{sr} (B) = 1$ and $K_1 (B) = 0.$ Then
$\Invo (\tilde{B})$ is dense in $\tilde{B}$, and it easily follows that
$\Invo (B)$ is dense in $B + 1.$ Next, note that
$[\tilde{B}, \tilde{B}] = B$ by Lemma~4.2.  Lemma~3.2 then implies that
$\exp(b) \in P(B+1)^-$ for all $b \in B$. Since every element of
$\Invo (B)$ is a product of exponentials $\exp(b)$ with $b \in B$,
it follows that $B + 1 = \Invo(\tilde{B})^- \subset P(B+1)^-$. \qed

We note that it is shown in~[\Leen] that every element of
$\Invo (K \otimes A)$ is a product of finitely many positive invertible
elements of ${(K \otimes A)}\tilde{}$. Clearly we may restrict these
elements to be in $(K \otimes A) + 1$. Thus, in the theorem above,
if $\oper{sr} (A) = 1$ and $K_1 (A) = 0$, then every invertible \el\  of
$(K \otimes A) + 1$ is in fact exactly a product of positive \el s.

\rem
This theorem, combined with Theorem 2.4, shows that if $A$ has
the \APFP , then so does $K \otimes A.$
The converse is false, even for simple \ca s.  It is
easy to construct, using the methods of~[\GdB] (as in Example~4.11), an
infinite-dimensional
simple unital \ca\  $A$, obtained as a direct limit with slow dimension
growth, such that $A$ has trivial $K_1$-group, stable rank $1$, and
real rank $1$. Theorem~4.14 implies that $K \otimes A$ has the
\APFP , while Theorem~4.9 implies that $A$ does not.
\qedd

\exa
We give an example of a unital \ca\  $A$  with the \APFP\  which
contains an ideal $I$ without the \APFP . This example shows that
Theorem~3.6 can't be extended to say that the \APFP\  passes to ideals,
even when the quotient also has the \APFP . Our algebra $A$ also does
not have real rank zero, showing that in general the \APFP\  does
not imply real rank zero.

Set $I_0 = K \otimes C(\tn)$. It follows from Kasparov's extension
theory~[\Ks] and the Universal Coefficient Theorem (Theorem 1.17
of~[\RS]) that there is an extension
$$ 0 \to I_0 \to A_0 \to \cn \to 0 $$
such that the induced map from $K_0 (\cn)$ to $K_1 (I_0)$ is an
isomorphism. (Note that $A_0$ will not be unital.) The long exact
sequence in K-theory then shows that $K_1 (A_0) = 0$, and Theorem
4.11 of~[\RfA] shows that $\oper{sr} (A_0) = 1$. (Note that
$\oper{sr} (I_0) = 1$.) So $K \otimes A_0$ has the \APFP\  by
Theorem 4.14. Its ideal $K \otimes I_0$ does not have the
\APFP , because $K_1 (K \otimes I_0) \neq 0$. Furthermore,
$K \otimes A_0$ does not have real rank zero, because it
contains the ideal $K \otimes I_0$ whose real rank is not zero.
(See~[\BP].)

The algebra $K \otimes A_0$ is not unital, but it can be used to
construct a unital example $B$ in the following two ways. Either
unitize it,
and then form the tensor product with an infinite dimensional UHF
algebra $F$,
or form a unital trivial extension by an infinite dimensional UHF
algebra.

In the first case, $B$ will have the \APFP\  by Theorem 4.6(1), and will
still contain the ideal $K \otimes I_0$ which has neither the
\APFP\  nor real rank zero. Since this ideal does not have real rank
zero, neither does $B.$ In the second case, $B$ will have the \APFP\  by
Theorem~3.6(1). It will now contain the ideal
$F \otimes K \otimes I_0 \cong F \otimes K \otimes K \otimes C(\tn).$
Since this ideal has nontrivial $K_1,$ it does not have the \APFP.
Cutting it down by a rank one projection in $K \otimes K,$ we obtain a
hereditary subalgebra of $B$ which is isomorphic to $F \otimes C(\tn).$
This algebra does not have real rank zero because the projections don't
distinguish the traces. It follows from Corollary 2.8 of [\BP] that
$B$ does not have real rank zero either.
\qedd

\address

\refs

\pref{C.S. Ballantine}{Products of positive definite matrices,
IV}{Linear Alg. Appl.}{3}{70}{79--114}

\pref{B. Blackadar, M. D\v{a}d\v{a}rlat, and M. R\o rdam}{The
real rank of inductive limit \ca s}{Math. Scand.}{69}{91}{211--216}

\pref{L.G. Brown and G.K. Pedersen}{C*-algebras of real rank zero}{J.
Funct. Anal.}{99}{91}{131--149}

\pref{H. Choda}{An extremal property of the polar decomposition in
von~Neumann \al s}{Proc. Japan Acad.}{46}{70}{341--344}

\pref{M.-D. Choi}{A simple \ca\  generated by two finite-order
unitaries}{Canadian J. Math.}{31}{79}{867-880}

\pref{J. Cuntz}{K-theoretic amenability for discrete groups}{J. reine
ang. Math.}{344}{83}{180-195}

\pref{J. Cuntz and G.K. Pedersen}{Equivalence and traces on
\ca s}{J. Funct. Anal.}{33}{79}{135--164}

\pref{P. de la Harpe and G. Skandalis}{D\'{e}terminant associ\'{e}
\`{a} une trace sur une alg\`{e}bre de Banach}{Ann. Inst.
Fourier}{43-1}{84}{241--260}

\pref{K.R. Goodearl}{Notes on a class of simple \ca s with real rank
zero}{Publ. Sec. Mat. Univ. Aut\'{o}noma Barcelona}{36}{92}{637--654}

\ppref{K.R. Goodearl}{Riesz decomposition in inductive limit
\ca s}{92}{}

\ppref{U. Haagerup}{Quasitraces on exact \ca s are
traces}{91}{(handwritten manuscript)}

\bref{P.R. Halmos}{A Hilbert Space Problem
 Book}{Springer}{New~York}{82}

\pref{G.G. Kasparov}{The operator K-functor and extensions of
\ca s}{Izv. Akad. Nauk SSSR Ser. Mat.}{44}{80}{571--636}

\pref{M. Khalkali, C. Laurie, B. Mathes, and H. Radjavi}{Approximation
by products of positive operators}{J. Operator Theory}{29}{93}{237--247}

\thesis{M. Leen}{Factorization in the Invertible
Group of a $C^*$-Algebra}{University of Oregon, Eugene}{94}

\pref{H. Lin}{Exponential rank of \ca s with real rank zero and
Brown-Pedersen's conjecture}{J. Funct. Anal.}{114}{93}{1--11}

\pref{G.J. Murphy}{The analytic rank of a C*-\al}{Proc. Amer. Math.
Soc.}{115}{92}{741--746}

\pref{C. Pearcy and D.M. Topping}{Sums of small numbers of
idempotents}{Michigan J.~Math.}{14}{67}{453--465}

\pref{N.C. Phillips}{Simple \ca s with the property weak (FU)}{Math.
Scand.}{69}{91}{127--151}

\pref{N.C. Phillips}{Exponential length and traces}{Proc. Royal Soc.
Edinburgh, Sec. A}{}{95(?)}{to appear}

\ppref{N.C. Phillips}{Factorization problems in the invertible
group of a homogeneous \ca}{93}{University of Oregon}

\ppref{T. Quinn}{Ideals in AF-algebras and approximation by products of
positive operators}{93}{Trinity College, Dublin}

\pref{H. Radjavi}{Products of Hermitian matrices and symmetries}{Proc.
Amer. Math. Soc.}{21}{69}{369--372; {\bf 26} (1970), 701}

\pref{M.A. Rieffel}{Dimension and stable rank in the K-theory of
C*-\al s}{Proc. London Math. Soc. {\rm (3)}}{46}{83}{301--333}

\pref{M.A. Rieffel}{The homotopy groups of the unitary groups of
noncommutative tori}{J. Operator Theory}{17}{87}{237--254}

\pref{M. R\o rdam}{On the structure of simple \ca s tensored with
a UHF algebra}{J. Funct. Anal.}{100}{91}{1--17}

\pref{M. R\o rdam}{On the structure of simple \ca s tensored with
a UHF algebra, II}{J. Funct. Anal.}{107}{92}{255--269}

\pref{J. Rosenberg and C. Schochet}{The K\"{u}nneth theorem and the
universal coefficient theorem for Kasparov's generalized K-functor}{Duke
Math. J.}{55}{87}{431--474}


\pref{P.Y. Wu}{Products of normal operators}{Canadian J.
Math.}{40}{88}{1322--1330}

\bye